\input amstex
\documentstyle{amsppt}
\vsize=8in \hsize 6.6 in 
\loadbold \topmatter
\title Direction control of bilinear systems. I. \endtitle
\author Abdon Choque,  Valeri Marenitch \endauthor

\address  Universidad Michoacana, Morelia, Mexico; Kalmar, Sweden\endaddress \email abdon\@ ifm.umich.mx;  valery.marenich\@ gmail.com \endemail
\keywords   \endkeywords \subjclass 93B05, 93B29, 20M20 \endsubjclass
\date{}\enddate
\abstract We consider the general bilinear control systems in three-dimensional Euclidean space which satisfy the $LARC$ and have an open control
set $U$; and give sufficient controllability conditions for corresponding projected systems on two-dimensional sphere.
\endabstract

\endtopmatter

\document

\head Introduction and Statement of Results \endhead

Bilinear control dynamical systems in a Euclidean space $x\in R^n$ are given by linear equations
$$
\dot x(t)=A(u)x(t), \qquad (\Sigma)\tag 1
$$
with the right-hand linear operator $A(u)$ depending linearly on the control parameters $u=(u^1,...,u^d)\in U \subset R^d$:
$$
A(u)=A + u^1 B_1 + ... + u^d B_d, \tag 2
$$
where $A$ and $B_k, k=1,...,d$ are some constant linear operators. Then $x(t)$ is the solution (or trajectory) of $\Sigma$ if
$$
\dot x(t)=A(u(t))x(t) \tag 3
$$
for some {\bf control} $u(t)$; i.e., the function $u(t): [0,T] \to U$. The positive orbit ${\Cal O}^+(x)$ of the point $x$ is the union of all
trajectories $x(t)$ issuing from it; i.e., such that $x(0)=x$. The bilinear system $\Sigma$ is called {\bf controllable} if for every non-zero
point $x$ its orbit coincide with the whole space without its origin zero:  ${\Cal O}^+(x)=R^n\backslash{0}$, which means that for two arbitrary
non-zero points $x,y$ there always exists some trajectory $x(t)$ of $\Sigma$ which connect them: $x=x(0), y=x(T)$. We may "project" the bilinear
system to the control dynamical system defined on the sphere $S^{n-1}$ of unit vectors in $R^n$ by considering an evolution of unit directions
$q(t)=x(t)/\|x(t)\|$ of the non-zero solutions $x(t)$ of $\Sigma$: for the unit vector $q$ we define the vector $V(q,u)$ to be equal to the
tangent to the sphere component of the vector $A(u)q$; i.e., $V(q,u)=A(u)q - (A(u)q,q)q$. Then
$$
\dot q(t)=V(q(t),u), \qquad (\Sigma^{pr}) \tag 4
$$
correctly defines the control dynamical system on the unit sphere which we call the {\bf projection} $\Sigma^{pr}$ of the system $\Sigma$. This
definition is justified since for an arbitrary (non-zero) trajectory $x(t)$ with the control $u(t)$ its central projection $x(t) \to x(t)/
\|x(t)\|$ to the unit sphere is the trajectory $q(t)$ with the same control of the projected system $\Sigma^{pr}$. If in addition some linear
operator $A(u_*)$ has all eigenvalues with negative real parts, while some other $A(u^*)$ has all eigenvalues with positive real parts (i.e., the
first one takes eventually an arbitrary big sphere inside the unit sphere, while the second - an arbitrary small sphere outside the unit sphere),
then the controllability of $\Sigma$ is equivalent to the controllability of its projection $\Sigma^{pr}$: in [CK] (see Corollary 12.2.6) the
authors the claims that the system $\Sigma$ is controllable if and only if

CK1) their Lyapunov spectrum contains zero in its interior, and

CK2) the projected system $\Sigma^{pr}$ of $\Sigma$ is controllable.

\medskip

Note, that from the linearity of $\Sigma$ it follows that $\Sigma^{pr}$ is invariant under the central symmetry $q\to -q$, and therefore, we may
further reduce it to the system on the projective space $RP^{n-1}=S^{n-1}/Z_2$. We do not do this, since for some symmetric and non-controllable
systems on the sphere their further reductions to the systems on projective space are controllable; i.e., this further reduction does not
preserves simultaneous controllability. It is, also, more convenient for us to work on the unit sphere in the Euclidean space.

\medskip

Next question is: when the projected system is controllable? The two-dimensional case $n=2$ was considered in details in [AS] and [BGRS] where the
authors show that the projected system on the circle $RP^1$ satisfying the $LARC$ is controllable if and only if for some parameter $u_0$ the
matrix $A + u_0 B$ has complex eigenvalues.\footnote{For the sake of completeness note that all one-dimensional ($n=1$) control systems $\dot
x=A(u)x$ are controllable in our sense iff $A(u)\not\equiv 0$ since all non-zero $x$ have the same direction - the set of unit direction $RP^0$ is
a one-point set.} Remind, that the $LARC$ - the "Lie Algebra Rank Condition" is satisfied for the control dynamical system $\Sigma$ at the point
$x$ if the evaluation at this point of the algebra Lie ${\Cal LA}(\Sigma)$ generated by all vector fields $A(u)x$ coincide with the whole $R^n$.
Easy to see that the systems $\Sigma$ and $\Sigma^{pr}$ satisfy or do not satisfy $LARC$ simultaneously, and that for $n=3$ the system
$\Sigma^{pr}$ on the unit two-dimensional sphere $S^2\subset R^3$ satisfies $LARC$ at $q$ if and only if there exists two linearly independent
vectors $V(q,u)$ and $V(q,v)$ at this point (our CC1 condition below).

\medskip

In general, the controllability of the bilinear system (1) corresponds to the transitivity of the Lie algebra generated by the right-hands $A(u)$.
Classification of such transitive actions is obtained in [B] (see also [BW], where the algorithm is given for the practical verification of this
property). Our approach is quite different: instead of checking the transitivity we study directly the chronological products of the control
system (1). The reason is the following. In [M1] we considered quantum systems with Hamiltonian depending on a parameter $u$. Using the Dirac
mechanism of quantization, it was suggested some new description of the Berry-phase geometry for such systems. The problem of measurement of such
systems naturally leads to the investigation of dynamical systems with the control $u$. These systems are defined on "level" sub-manifolds $M$ of
the Hamiltonian in the space of parameters; and are not linear systems. We verify the controllability of such systems with the help of the "closed
orbit controllability criterium", see [M2], which works for arbitrary non-linear control systems on closed manifolds. For instance, structurally
stable dynamical systems on any closed surface $M^2$ have decompositions into invariant cells with the same dynamics, as we describe below in this
paper; and their controllability follows from the same incidence relations. For bilinear systems these relations are extremely simple - they are
given by linear inequalities, which is presented in this paper.

More precisely, the main goal of this paper is to give sufficient conditions for controllability of bilinear systems satisfying the $LARC$ (CC1
condition) in three-dimensional case and have an open control set $U$ (CC2 condition). We do not consider yet the case $n=3$ completely in all
details, but restrict our attention to the "stable" or "general" systems (see [PDM]) as follows: we choose some finite subset ${\Cal U}^N$ of
parameters $u_k, k=1,...,N$ and assume (our stability or generality conditions CC3, CC4; see below) that the systems of linear equations
$\Sigma_k$ for these parameters have different eigenvalues and their eigenvectors $E^k_i$ are different. Then we define - by a number of linear
conditions on these eigenvectors - the open subsets $CO^2(E^k_i)$, which consist of geodesic triangles in the unit sphere. Finely, we prove that
the projected system $\Sigma^{pr}$ is controllable if some simple incidence relations are satisfied, see the Theorems~A-C below.

\medskip
Our arguments rely on the "Closed Orbit Controllability Criterium" from [M2], where it is proved that the control dynamical system on a compact
manifold satisfying the $LARC$ and with an open set of control vectors (see definitions in [M2]) is controllable if and only if through an
arbitrary point goes some closed trajectory of the system. To make this paper self-content we reprove the two-dimensional version of this
criterium here. Then we find sufficient controllability conditions such that compositions of shifts $P^t_k$ along trajectories of some finite
number $A(u_k)$ provide us with the collection of closed trajectories going through an arbitrary $q$.

The authors sincerely thank Professor D.~Elliott for valuable comments.

\medskip

\head 1. Bilinear control systems in $R^3$ \endhead

As we mentioned above, the question if the given bilinear system $\Sigma$ in $R^3$ is controllable or not can be reduced to the same question
about controllability of the projected system $\Sigma^{pr}$ on the sphere $S^2$. In this section we consider first the more general question of
the controllability of systems $\Sigma_M$ on closed (connected, compact without boundary) surfaces $M^2$ - we reprove for $n=2$ the closed orbit
controllability criterium from [M2]. Then, using this criterium, we give sufficient conditions for the controllability of the projected systems
$\Sigma^{pr}$ on the two-sphere $S^2$. These conditions we formulate in terms of finite control subsystems of $\Sigma^{pr}$ - for finite number of
parameters $u_k, k=1,...,N$ we consider systems
$$
\dot q(t)=V(q(t),u_k) \qquad (\Sigma^{pr}_k), \tag 5
$$
and find conditions on the set of $\Sigma^{pr}_k$ such that compositions of shifts $P^t_k$ along trajectories of these systems (i.e.,
chronological products, see [M2]) give us closed trajectories going through all points $q$ of the sphere. This is equivalent to say that every two
points $q$ and $p$ may be connected by some trajectory of the system with piece-wise constant control functions with values $u_k, k=1,...,N$.

\medskip

\head 2. Controllability conditions \endhead

The first condition we require is that the collection of vector fields $V_k(q)=V(q,u_k), k=1,...,N$ satisfies the $LARC$ which for $S^2$ or any
other surface $M^2$ is the following:

\medskip

{\bf CC1. In an arbitrary point $q\in M^2$ the set of all control vectors $\{V(q)\}$ contains two linearly independent vectors.}

\medskip

For the control dynamical system on the surface $M^2$ the set $V(q)=\{V(q,u)| u\in U\}$ of all control vectors at the point $q$ is the subset of
the tangent plane $T_qM^2$ to $M^2$ at the point $q$. The union $V$ of all $V(q)$ is the subset of the union of all planes $T_qM^2$ - the tangent
bundle $TM^2$. The second condition we require is

\medskip

{\bf CC2. The set $V$ is an open subset in the tangent bundle $TM^2$.}

\medskip

For projections $\Sigma^{pr}$ of bilinear systems with an open control set $U$ this second condition CC2 follows from the first one; i.e., from
the $LARC$.

\medskip

\proclaim{Lemma~1}(Closed orbit controllability criterium for $n=2$) Let $\Sigma_M$ be some control dynamical system on the closed connected
surface $M^2$:
$$
\dot q(t)=V(q(t),u). \tag 6
$$
If the system $\Sigma_M$ satisfies CC1-CC2 above, then it is controllable if and only if through an arbitrary point $q$ goes some closed and
non-trivial orbit of the system $\Sigma_M$.
\endproclaim

\demo{Proof} Let $B$ be some coordinate neighborhood of the point $q$ with local coordinate functions $x^i: B \to R, i=1,2$ such that $x^i(q)=0$.
Then the tangent plane to $M^2$ at the point $q'\in B$ is generated by two coordinate vectors $e_i(q'), i=1,2$ - evaluations at this point $q'$ of
the coordinate vector fields $e_i=\partial / \partial x_i$, while for the control vectors holds:
$$
V(q',u)=V^1(q',u)e_1(q') + V^2(q',u)e_2(q') \tag 7
$$
for some functions $V^i(q',u), i=1,2$. If $q(t), 0\leq t\leq T$ is the closed trajectory of $\Sigma_M$ with the control function $u(t)$ going
through the point $q$, $q=q(0)=q(T)$, then
$$
{{d}\over{dt}}x^i(q(t))=V^i(q(t),u(t)).\tag 8
$$
Because CC2 is satisfied, for sufficiently small positive $\epsilon$ all vectors $V'$ from $T_{q(t)}M$ for $0<t<\epsilon$ which are
$\epsilon$-close to $V(q(t),u(t))$ belong to the control set $V(q(t))$. Therefore, for an arbitrary small vector $v=(v^1,v^2)$ with
$(v^1)^2+(v^2)^2< \epsilon^2$ the solution of the equation
$$
{{d}\over{dt}}x_v^i(t)=V^i(q(t),u(t)) + v^i \tag 9
$$
also gives us the trajectory $q_v(t)=(x_v^1(t),x_v^2(t))$ of the control system $\Sigma_M$. Since, for a given $0<t<\epsilon$ the set of all
points $q_v(t)$ is an open set which contains $q(t)$, we see that the set of points ${\Cal O}^+(q)$ which are reachable from $q$ contains some
open neighborhood of (every) $q(t)$ for  $0<t<\epsilon$. By standard theorems on ordinary differential equations, the image of every such open
neighborhood under all shifts $P^\tau, 0<\tau<T-t$ along trajectories of the system $\Sigma_M$ with the control function $u(t)$ is an
open-neighborhood of the trajectory $q(t), 0<t<T$. Therefore, the positive orbit ${\Cal O}^+(q)$ of the point $q$ contains together with arbitrary
closed trajectory going through $q$ some its open neighborhood - which we denote by $B^+(q)$. The negative orbit ${\Cal O}^-(q)$ of $q$ is the set
of points from which we can reach the point $q$ by trajectories of our system $\Sigma_M$. It is by definition the positive orbit of the
"minus"-$\Sigma_M$ control system, i.e., the system $-\Sigma_M$ on $M$ with the control set $-V(q)$. Easy to see that these two systems satisfy or
do not satisfy our conditions CC1-CC2 simultaneously. Therefore, we may repeat our arguments and prove that ${\Cal O}^-(q)$ contains with an
arbitrary point $q(t)$ for negative $t<0$ some its open neighborhood $B^-(q(t))$. Because for our closed orbit it holds $q(t)=q(T+t)$ we see that
for arbitrary $0\leq \tau\leq T+t$ this neighborhood is an image of some open neighborhood of $q(\tau)$ under the shift $P^{T-\tau+t}$ along
trajectories of $\Sigma_M$. Which implies that there exists some open neighborhood $B^-(q)$ of the closed orbit such that from its every point we
can reach $q$. Then by definition the system $\Sigma_M$ is controllable on the intersection $B_q=B^-(q)\cap B^+(q)$, which is again an open
neighborhood of the closed orbit $q(t), 0\leq t \leq T$.

Now controllability of $\Sigma_M$ on the connected surface $M$ follows from its compactness, see [M2]. For every $q\in M$ there exists some closed
orbit going through this point. Every such orbit has an open neighborhood $B_{q}$ on which the system $\Sigma_M$ is controllable. The union of all
$B_{q}$ covers $M^2$. Since $M^2$ is compact there exists some finite sub-covering
$$
M^2=\cup_{i=1,...,c} B_{q_i}. \tag 10
$$
Because $M^2$ is connected for two arbitrary points $q$ and $p$ we may find a finite sequence $B_{q_{i_1}},..., B_{q_{i_k}}$ from the finite
covering (10) such that $q\in B_{q_{i_1}}$, $p\in B_{q_{i_k}}$ and $B_{q_{i_s}} \cap B_{q_{i_{s+1}}} \not= \emptyset$. Take $q_0=q$; some $q_s \in
B_{q_{i_s}} \cap B_{q_{i_{s+1}}}$ and $q_{k+1}=p$. Since $\Sigma_M$ is controllable on each $B_{q_{i_s}}$ there exists a trajectory from $q_s$ to
$q_{s+1}$ for all $s=0,..., k$. Taking the composition of these trajectories we obtain the trajectory of $\Sigma_M$ going from $q$ to $p$. Since
$q,p$ where arbitrary this proves that $\Sigma_M$ is controllable on $M^2$.

To complete the proof of our controllability criterium (... "and only if" part) we note that if the system is controllable than for arbitrary $q$
there exists the trajectory $\gamma$ of the system going from $q$ to the point $q(-\epsilon)$ for some other non-trivial trajectory $q(t),
-\epsilon\leq t\leq \epsilon$ of the system $\Sigma_M$ going through $q=q(0)$. Then the composition of $\gamma$ and $q(t), -\epsilon < t \leq 0$
gives us the non-trivial closed trajectory through the point $q$. The Lemma is proved.
\enddemo

\medskip

Now we go back to the projected bilinear systems.

For every fixed $u_k$ the system $\Sigma^{pr}_k$ is the projection of the system of linear differential equations in $R^3$:
$$
\dot x(t)=A(u_k)x(t) \qquad (\Sigma_k). \tag 11
$$
Denote by $\lambda^k_1(u),\lambda^k_2(u),\lambda^k_3(u)$ eigenvalues of $A(u_k)$ and by $E^k_1(u),E^k_2(u),E^k_3(u)$ the corresponding unit
eigenvectors.

Here we will split our consideration into two cases: first is when all eigenvalues are real, and the second one - when one $\lambda^k_R$ is real
and other two are complex conjugate numbers $\lambda^k_C$ and $\overline{\lambda^k_C}$

\medskip

\head 3. Real eigenvalues. Jordan cells decomposition and dynamics \endhead

Here we further assume

{\bf CC3. The eigenvalues $\lambda^k_i,i=1,2,3$ are different:}

$$
\lambda^i_1<\lambda^i_2<\lambda^i_3,  \tag 12
$$
which means that the linear operator $A(u_k)$ is stable - similar to all other operators which are close enough to it. Then in the basis
$E^k_i(u),i=1,2,3$ the solutions of the system (11) are given by
$$
x_k(t)=(x_k^1(t)=x_k^1(0) e^{\lambda^k_1 t}, x_k^2(t)=x_k^2(0) e^{\lambda^k_2 t}, x_k^3(t)=x_k^3(0) e^{\lambda^k_3 t}). \tag 13
$$
We see that closer of domains in $R^3$ where coordinates have the same sign are invariant sets of this system. Correspondingly, triangles on the
unit sphere with vertices $\pm E^k_i,i=1,2,3$ are invariant subsets of the projected system $\Sigma_k^{pr}$. Sides of these triangles (also
invariant subsets) are intervals of big circles (geodesics) in the unit sphere $S^2$. We denote by $C^k_i$ the circle around $E^k_i$ consisting of
such intervals; e.g., the circle $C^k_3$ consists of the intervals from $E^k_1$ to $E^k_2$, then from $E^k_2$ to $-E^k_1$, then from $-E^k_1$ to
$-E^k_2$ and from $-E^k_2$ to $E^k_1$. Inside each invariant triangle the position of the point $q$ may be given with the help of coordinates
$x_k$. For definiteness we consider the triangle $\Delta^k=\Delta^k(+++)$ with vertices $E^k_1,E^k_2,E^k_3$, where all $x^i_k$ are non-negative.
For an arbitrary point $q$ inside this triangle its coordinates $(x_k^1(q),x_k^2(q),x_k^3(q))$ in the coordinate system $x_k$ provide homogeneous
coordinates $(x_k^1(q):x_k^2(q):x_k^3(q))$ for the unit direction $q$ in the projective sphere $RP^{2}$; it holds
$(x_k^1(q))^2+(x_k^2(q))^2+(x_k^3(q))^2=1$. The vector field $V_k(q)=V(q,u_k)$ of the dynamical system $\Sigma^{pr}_k$ due to (12) is stable (see
[PDM]) and has six zeroes: two sources at $\pm E^k_1$, two saddles at $\pm E^k_2$ and two sinks at $\pm E^k_3$. Circles $C^k_3$ consist of saddle
separatrix incoming to saddles $\pm E_2^k$, and circles $C^k_1$ - of outgoing separatrix from saddles to sinks $\pm E^k_3$.

If we denote by $P^t_k$ the shift along trajectories of $\Sigma^{pr}_k$ - the one-parameter group of diffeomorphisms of $S^2$ generated by the
vector field $V_k$ and given by solutions of (11) - then from (13) (considering mutual ratios of coordinates $x_k^i(t)$ as $t\to\infty$) we easily
conclude the following description of the dynamic of such shifts inside our invariant triangle $\Delta_k$.

\medskip

\proclaim{Lemma~2}

1.1. For an arbitrary point $q$ in the interior of $\Delta_k$ and arbitrary small $\epsilon$-neighborhood $B_1^k(\epsilon)$ of $E_1^k$ there
exists a trajectory of $\Sigma_k^{pr}$ starting at some point $q(\epsilon)$ in $B_1^k(\epsilon)$, going through $q$ and converging to $E_3^k$.

1.2. For an arbitrary curve $q(t), -\epsilon<t<\epsilon$ which is transversal to the outgoing separatrix $E_2^kE_3^k$ and intersect it at the
point $q=q(0)$, and another curve $\bar q(\bar t), -\bar\epsilon<\bar t <\bar\epsilon$ which is transversal and intersect the incoming separatrix
$E^k_1E^k_2$ at the point $\bar q=\bar q(0)$, such that both $q(t)$ and $\bar q(\bar t)$ belong to the interior of the triangle $\Delta_k$ for
positive $t,\bar t$, for small enough positive $t<\epsilon'$ there exists a trajectory of $\Sigma_k^{pr}$ starting at some point $\bar q(\bar t)$
for some positive $\bar t=\bar t(t)$ and going through the point $q(t)$.
\endproclaim

\medskip

Note that the claims of the Lemma~2 are valid also for non-stable dynamical systems when (12) is wrong, but equal eigenvalues correspond to
$A(u_k)$ given by the Jordan form matrix on the subspace generated by corresponding eigenvectors. In other two cases we have or $V_k$ identically
zero - when $\Sigma_k$ is diagonal and has equal eigenvalues; or $V_k$ equals zero on some side of the triangle $\Delta_k$ corresponding to the
two-dimensional subspace, where the action of $A(u_k)$ is diagonal with equal eigenvalues.

\medskip

\head 4. Complete subsystems. Real eigenvalues\endhead

Let ${\Cal U}^N=\{u_k|k=1,...,N\}$ be some finite subset of the set $U$ of control parameters, and $\Sigma^{pr}_k$ as above denotes the projection
(11) of $\Sigma_k$. We call ${\Cal U}^N$ {\bf complete} if every two points $q,p$ on the unit sphere can be connected by some trajectory of
$\Sigma^{pr}$ with the piece-wise control function $u(t):[0,T]\to {\Cal U}^N$; i.e., taking values in the finite subset ${\Cal U}^N$. Now we give
some sufficient conditions which guarantee that ${\Cal U}^N$ is complete. First we assume in addition to CC3 that our subsystem is general:

\medskip
{\bf CC4. All eigenvectors $E^k_i$ are different, and $E^k_i$ does not belong to any plane $\Pi^j_{s_1}$ generated by $E^j_{s_2}$ and $E^j_{s_3}$
for all different $s_1,s_2,s_3=1,2,3$.}
\medskip

Next, for an arbitrary point $q$ we generate two sets: first is the set $CO^1(q)$ of saddle points $E^j_2$, and the second set $CO^2(q)$ - some
open subset of the unit sphere, which is the union of (geodesic) triangles in this sphere with sides being big circles. Finely, our completeness
criterium is that for some index $s$ both $E^s_1$ and $-E^s_1$ belong to $CO^2(E^s_3)$ and $CO^2(-E^s_3)$, see the Theorem~A below. This criterium
is given in terms of a number of simple inequalities, which describe mutual positions of systems of eigenvectors $\{E^k_i,i=1,2,3\}$. The
construction of the sets $CO^1$ and $CO^2$ proceeds as follows.

\medskip

{\bf Construction of $CO^1$ and $CO^2$} is by repetition of steps of two types.

{\bf First type step.}

Take some point $q$ in the interior of $\Delta_k$ and assume that some circle $C^j_3$ "cut" this point from the vertex $E^k_3$ of the triangle
$\Delta_k$. Then by the Lemma~2 above the trajectory $q(t)=P^t_k(q)$ of $\Sigma^{pr}_k$ issuing from $q$ converge to $E^k_3$, and therefore,
intersect the circle $C^j_3$ at some point $\bar q=P^{\bar t}_k(q)$. Then for $t<\bar t$ and close to $\bar t$ the points $q(t)$ lie in some, say
$\Delta_j(++-)$ - the triangle with vertices $E_1^j, E^j_2,-E^j_3$, while for $t>\bar t$ and close to $\bar t$ the points $q(t)$ lie in the
triangle $\Delta_j(+++)$ - the triangle with vertices $E_1^j, E^j_2,E^j_3$. Applying the Lemma~2 again, now for the triangles $\Delta_j(++-)$ and
$\Delta_j(+++)$  we see that with the help of shifts along trajectories of $\Sigma^{pr}_j$ we may approach arbitrarily close an arbitrary point of
the semi-circle consisting of outgoing separatrix from $E^j_2$  - the union of intervals $-E^j_3E^j_2 \cup E^j_2E^j_3$, which we denote
$Sep^j_2(+)$. Since our trajectories under the flow $P^{t'}_j$ are coming from a neighborhood of the point $\bar q$, they converge to the
separatrix only from one side, the points which are reachable from $\bar q(t)$ always belong to the same half-sphere relative to $C^j_3$ as the
point $\bar q$. Note in addition, that for the given point $q$ in the interior of $\Delta_k$ and separated from $E^k_3$ by $C^j_3$ we may find
some open neighborhood $B(q)$ such that the same is true for all $q'\in B(q)$. Thus, from the Lemma~2 we conclude the following.

\medskip
\proclaim{Lemma~3} If the point $q$ in the interior of $\Delta_k$ is separated from $E^k_3$ by $C^j_3$, then for all $q'$ in some small
neighborhood $B(q)$ of $q$, and an arbitrary point $p$ on $Sep^j_2(+)$ which is not end point of it, i.e., $p\not= \pm E^j_3$ there exist an open
neighborhood $B(p)$ such that all points in $B(p)$ are reachable from $q'$ if they belong to the same half-sphere relative to $C^j_3$ as the point
$\bar q$.
\endproclaim
\medskip

Let us write down our conditions in an analytic form. If $S(+)$ and $S(-)$ denote the two half-spheres in which the big circle $C$ divides the
unit sphere, then the point $q$ belongs to one of them depending on the sign of the scalar product of $q$ with the vector normal to the plane in
which this circle lies. Or, if $C$ lies in the plane $\Pi$ generated by $E_1$ and $E_2$ - the vector normal to $\Pi$ is the vector product
$[E_1,E_2]$, and the points $q$ and $p$ belong to different half-spheres if and only if the corresponding products have different signs:
$$
(q, [E_1, E_2]) (p, [E_1, E_2]) < 0,
$$
or in a more symmetric form:
$$
det(q, E_1, E_2) det(p, E_1, E_2) < 0,
$$
where by $det(X,Y,Z)$ we denote the mixed product of three vectors $X,Y,Z$. Define the function: $s(q,C,p)=1$ if (14) is satisfied, and
$s(q,C,p)=-1$ if not. Then the points $q$ and $E^k_3$ are separated by $C^j_3$ if
$$
s(q,C^j_3,E^k_3) = 1. \tag 14
$$

The triangle $\Delta_k$ is the intersection of three half-spheres out of six $S^k_{i}(\pm), i=1,2,3$ in which the big circles $C^k_i, i=1,2,3$ -
its sides - divide the unit sphere. So, the condition that the point $q$ belong to the interior of $\Delta_k$ is that the vector $q$ belongs to
the same half-sphere relative to $C^k_i$ as the vertex $E^k_i$ (i.e., which does not belong to this circle):
$$
(q, [E^k_1, E^k_2])(E^k_3, [E^k_1, E^k_2]) > 0, (q, [E^k_2, E^k_3]) (E^k_1, [E^k_2, E^k_3]) > 0, (q, [E^k_3, E^k_1])(E^k_2, [E^k_3, E^k_1]) > 0,
$$
or
$$
det (q, E^k_1, E^k_2) det (E^k_3, E^k_1, E^k_2) < 0, det (q, E^k_2, E^k_3) det (E^k_1, E^k_2, E^k_3) < 0, det (q, E^k_3, E^k_1) det (E^k_2, E^k_3,
E^k_1) < 0. \tag 15
$$
Again, we define the function $d(q,\Delta)=1$ if the point $q$ belongs to the interior of the triangle $\Delta=\Delta(E_1,E_2,E_3)$ with vertices
$E_1,E_2,E_3$; and $d(q,\Delta)=-1$ - if not. From the above we see that $d(q,\Delta(E_1,E_2,E_3))=1$ if and only if
$$
s(q,C_i,E_i)=-1 \qquad \text{ for all } \qquad i=1,2,3, \tag 16
$$
where $C_i$ denotes the big circle in the sphere going through the points $E_j,j\not= i$.

If the incidence condition (15) holds for some $k$, and for this $k$ the separation condition (14) holds for some other $j$, then we set:
$$
E^j_2 \in CO^1(q). \tag 17
$$

{\bf Second type step.}

Let the point $E^j_2\in CO^1(q)$ belongs to some triangle $\Delta_l$, say $\Delta_l(+++)$, such that the semi-circle $Sep^j_2(+)$ cut from this
triangle some vertex $E^l_s$; i.e., such that end points $\pm E^j_3$ of $Sep^j_2(+)$ lie outside $\Delta_l$. Then, depending on $s=1,2,3$ the
shifts $P^{t"}_l$ along trajectories of $\Sigma^{pr}_l$ of the points $p$ from the chord $\Gamma^j_l=Sep^j_2(+) \cap \Delta_l$ of the semi-circle
$Sep^j_2(+)$, which is inside the triangle $\Delta_l$, cover some open subregion $\Omega^j_{l}$, which we include in $CO^2(q)$ (and define it
precisely in a moment below). Since such points $p$ belong to the interior of $\Delta_l$ every trajectory $P^{t"}_l(p), 0\leq t"\leq T"$ can be
continued over some open interval $-\epsilon"<t"<0$ such that the extended trajectory contains points from both sides of the chord $\Gamma^j_l$.
Because all points $p$ from $\Gamma^j_l$ are interior points of the separatrix $Sep^j_2(+)$ (end points $\pm E^j_3$ of $Sep^j_2(+)$ lie outside
$\Delta_l$) we conclude from the Lemma~3 above that all points of the open domain $\Omega^j_l$ are reachable by trajectories from $q$ and from all
points $q'$ of some its small neighborhood $B(q)$.

In addition, when $\Gamma^j_l$ cut the vertex $E^l_s$ for $s\not= 3$ - it has some end point on the incoming separatrix $E^l_1E^l_2$. Therefore,
arguing as in the Lemma~3, we also include the limit point of this separatrix $E^l_2$ in $CO^1$ under this step of the second type.

Write down the analytic form of our arguments. The point $E^j_2$ belong to $\Delta_l(E^l_1E^l_2E^l_3)$ if
$$
d(E^j_2, \Delta_l(E^l_1E^l_2E^l_3)) = 1. \tag 18
$$
The endpoints of the semi-circle $Sep^j_2(+)$ do not belong to this triangle if
$$
d(\pm E^j_3, \Delta_l(E^l_1E^l_2E^l_3)) = -1. \tag 19
$$
Then the chord $\Gamma^j_l$ cut from this triangle the vertex $E^l_1$ if
$$
s(E^l_1,C^j_3,E^l_2)=1 \qquad \text{ and } \qquad s(E^l_1,C^j_3,E^l_3)=1. \tag 20
$$
In this case the union of all trajectories of $\Sigma^{pr}_l$ issuing from points of $\Gamma^j_l$ covers the part of $\Delta_l$ lying in the same
half-sphere relative to $C^j_3$ as $E^l_3$:
$$
\Omega^j_l=\{p| d(p,\Delta_l)=1, s(p,C^j_3,E^l_3)=-1 \}. \tag 21
$$
If the chord $\Gamma^j_l$ cut from this triangle the vertex $E^l_2$; i.e.,
$$
s(E^l_1,C^j_3,E^l_2)=1 \qquad \text{ and } \qquad s(E^l_2,C^j_3,E^l_3)=1,  \tag 22
$$
then the union of all trajectories of $\Sigma^{pr}_l$ issuing from points of $\Gamma^j_l$ contain the part of $\Delta_l$ lying in the same
half-sphere relative to $C^j_3$ as $E^l_2$ - small geodesic triangle with the vertex $E^l_2$ and the opposite side $\Gamma^j_l$. Indeed, an
arbitrary point $p$ inside this triangle by the Lemma~2 lies on some trajectory $p(t)$ coming from some small neighborhood of $E^l_1$ - outside
this triangle, and therefore intersecting $\Gamma^j_l$ at some point. Hence, the points of
$$
\Omega^j_l=\{p| d(p,\Delta_l)=1, s(p,C^j_3,E^l_2)=-1 \}. \tag 23
$$
are reachable from $\Gamma^j_l$ at this case (22). If the chord $\Gamma^j_l$ cut from this triangle the vertex $E^l_3$; i.e.,
$$
s(E^l_1,C^j_3,E^l_3)=1 \qquad \text{ and } \qquad s(E^l_2,C^j_3,E^l_3)=1,  \tag 24
$$
then similarly all the points from the small triangle with the vertex $E^l_3$ and the opposite to it side $\Gamma^j_l$ are reachable from the
points of this side, and
$$
\Omega^j_l=\{p| d(p,\Delta_l)=1, s(p,C^j_3,E^l_3)=-1 \}. \tag 25
$$
After this step we add to $CO^2(q)$ the open region $\Omega^j_l$ and to the list $CO^1(q)$ the vertex $E^l_2$ if (20) or (22) holds. The result of
the step of the second type we formulate as follows.

\medskip
\proclaim{Lemma~4} Let $E^j_2\in CO^1(q)$ belongs to some triangle $\Delta_l$, and the semi-circle separatrix $Sep^j_2(+)$ cut from this triangle
some vertex $E^l_s$; i.e., such that end points $\pm E^j_3$ of $Sep^j_2(+)$ lie outside $\Delta_l$. Then, depending on $s=1,2,3$ the open set
$\Omega^j_l$ defined by (21) if $s=1$, (23) if $s=2$ or (25) if $s=3$ is contained in the positive orbit ${\Cal O}^+(q')$ for all points $q'$ from
some small neighborhood $B(q)$ of $q$. If $s=1,2$ then in addition the conclusion of the Lemma~3 is true for all interior points $p$ of the
separatrix $Sep^l_2(+)$.
\endproclaim
\medskip

Now we give our sufficient controllability criterium for the subsystem of $\Sigma^{pr}_k, k=1,...,N$ where all eigenvalues $\lambda^k_i$ are real.

\medskip
\proclaim{Theorem~A} Let $\Sigma^{pr}$ be some projected system which satisfies the $LARC$ (CC1), the openness condition CC2) and generality
conditions CC3, CC4. Assume that for some finite subset ${\Cal U}^N=\{u_k|k=1,...,N\}$ of $U$ and some index $s$ after some finite number of
steps, each of which is of the first or second type as above, it holds
$$
E^s_1 \text{ and } -E^s_1 \text{ belong to } CO^2(E^s_3) \text{ and } CO^2(-E^s_3). \tag 26
$$
Then the system $\Sigma^{pr}$ is controllable.
\endproclaim

\demo{Proof} Take an arbitrary point $q$. As in the lemma~1 we see that, from CC1 and CC2 it follows that both the positive and negative orbits of
$q$ contain some open sets - neighborhoods $B^+(q(\epsilon))$ and $B^-(q(-\epsilon))$ of some trajectory $q(t)$ going through $q$. Thus, taking
partition of the unit sphere into the set of triangles $\Delta_s(\pm \pm \pm)$ we see that there is some of these triangles, denote it by
$\Delta_s^+$ with the vertex $E^s_3$ or $-E^s_3$ - say $E^s_3$ for definiteness; which intersect the positive orbit of $q$, say at some point
$q'$. Therefore, by the lemma~2 there exists some trajectory from $q'$ to some $p'$ in an arbitrary small neighborhood of this vertex $E_3^s$. If
we take this small neighborhood to be $B(E^s_3)$ as in the Lemmas~3 and 4 above, then we find a trajectory from $q$ to some point $p'$ from which
we can reach all the points of $CO^2(E^s_3)$ and $E^s_1$ and $-E^s_1$ in particular. Since $CO^2(E^s_3)$ is an open set, we can reach all the
points from some small neighborhoods of $E^s_1$ and $-E^s_1$, which we denote $B'(E^s_1)$ and $B'(-E^s_1)$.

Since the negative orbit of $q$ also contains some open subset, it has non-empty intersection, say the point $q"$, with the interior of some of
the triangles $\Delta_s^-$ with the vertex $E^s_1$ or $-E_1^s$, say $E^s_1$ for definiteness. Again, by the Lemma~2 the point $q"$ belongs to some
trajectory issuing from some point $p"$ which we may choose in an arbitrary small neighborhood of this vertex. Finely, the combination of these
orbits: from $q$ to $q'$, then to $p'$, then to $p"$, then to $q"$ and back to $q$ gives us the closed orbit going through $q$. Since we choose
the point $q$ in an arbitrary way, we have proved that the control dynamical system $\Sigma^{k}, k=1,...,N$ satisfies the closed orbit
controllability criterium and is controllable by the Lemma~1.

The Theorem~A is proved.
\enddemo

\medskip

\head 5. Complete subsystems. Complex eigenvalues \endhead

When the system $\Sigma_j$ has one real eigenvalue $\lambda^j_1=\lambda^j_R$ and the pair of complex-conjugate eigenvalues
$\lambda^j_2=Re(\lambda^j_C)+\sqrt{-1}Im(\lambda^j_C)$ and $\lambda^j_3=Re(\lambda^j_C) - \sqrt{-1}Im(\lambda^j_C)$ then its solutions are given
by $x_j(t)=(x_j^1(t),x_j^2(t),x_j^3(t))$ where
$$
x^1_j(t)= x_j^1(0) e^{\lambda^j_R t}
$$
and
$$
x^2_j(t)=x_j^2(0)cos(Im(\lambda^j_C)t)+ x_j^3(0)sin(Im(\lambda^j_C)t))e^{Re(\lambda^j_C)t},
$$
$$
x^3_j(t)= -x_j^2(0)sin(Im(\lambda^j_C)t)+ x_j^3(0)cos(Im(\lambda^j_C)t))e^{Re(\lambda^j_C)t}.
$$
Here we have only two two-dimensional invariant subsets - half-spaces $R^3_+=\{ x_j^1\geq 0\}$ and $R^3_-=\{ x_j^1\leq 0\}$. Correspondingly, the
projected system $\Sigma^{pr}_j$ has two invariant cells which are half-spheres. They intersect along the circle $C^j$ which is also an invariant
set of $\Sigma^{pr}_j$. Another invariant set - two zeros $\pm E^j_1$ of $V_j$. If
$$
\lambda^j_R > Re(\lambda^j_C), \tag 27
$$
then zeroes $\pm E_1^j$ are sources of the dynamical system $\Sigma^{pr}_j$, while $C^j$ is an attracting cycle, and if
$$
\lambda^j_R < Re(\lambda^j_C), \tag 28
$$
then zeroes $\pm E_1^j$ are sinks of the dynamical system $\Sigma^{pr}_j$, while $C^j$ is a repulsive cycle. The analog of the Lemma~2 describing
the dynamic of $\Sigma^{pr}_k$ is the following.

\medskip
\proclaim{Lemma~5}

1. If $E^j_1$ is the source, and $C^j$ is an attracting cycle, then for an arbitrary point $q$ in the interior of $\Delta_j$ and arbitrary small
$\epsilon$-neighborhood $B(E^j_1,\epsilon)$ of $E_1^j$ there exists a trajectory $q(t)$ of $\Sigma_j^{pr}$ starting at some point $q(\epsilon)$ in
$B(E^j_1,\epsilon)$, going through $q$ and converging to the big circle $C^j$: for an arbitrary point $p$ of this circle and arbitrary $q$ there
exists a sequence $t_n(q,p)\to \infty$ such that $q(t_n)$ converge to the point $p$.

2. Let $E^j_1$ be the sink, $C^j$ - the repulsive cycle, and $p(\bar t)$ some curve intersecting transversally $C^j$ at the point $p=p(0)$. Then
for an arbitrary point $q$ in the interior of $\Delta_j$ and arbitrary small $\epsilon$-neighborhood $B(p,\epsilon)$ of $p$ there exists a
trajectory $q_{\bar t}(t)$ of $\Sigma_j^{pr}$ starting at the point $p(\bar t)$ in $B(p,\epsilon)$, going through $q$ and converging to the sink
$E^j_1$.
\endproclaim

\medskip

Using arguments from the Theorem~A above we see that the $LARC$ (CC1) and the openness (CC2) conditions with the help of the Lemma~5 immediately
imply that the positive orbit of the source $E^j_1$ (under (27)) contains some open neighborhood of the half-sphere $\Delta^j(+)$ (and ${\Cal
O}^+(-E^j_1)\supset \Delta^j(-)$ correspondingly); while if $E^j_1$ is a sink (when (28) holds) even more is true: the positive orbit of any small
curve $p(\bar t)$ intersecting $C^j$ transversally at this point $p$ contains the whole unit sphere.

Having at least one such $\Sigma^{pr}_j$ in our finite subsystem simplifies our controllability considerations: Let ${\Cal U}_R^N$ denotes the
parameters $u_k$ in ${\Cal U}^N$ such that $\Sigma_k$ has only real eigenvalues, and $CO^1(q)$, $CO^2(q)$ constructed as above for the subsystem
${\Cal U}_R^N$. Then repetition of the arguments in the proof of the Theorem~A give us the following.

\medskip
\proclaim{Theorem~B} Let $\Sigma^{pr}$ be some projected system which satisfies the $LARC$ (CC1), the openness condition CC2) and generality
conditions CC3, CC4; and ${\Cal U}^N=\{u_k|k=1,...,N\}$ some finite set of parameters such that for one of them $j$ the system $\Sigma_j$ has
complex eigenvalues. Assume that for the subset ${\Cal U}_R^N=\{u_k|k=1,...,N\}$ of parameters such that all eigenvalues of $\Sigma_k$ are real,
and some index $s$ after some finite number of steps, each of which is of the first or second type as above, it holds:
$$
E^j_1 \text{ and } -E^j_1 \text{ belong to } CO^2(E^s_3) \text{ and } CO^2(-E^s_3),  \tag 29
$$
if (27) is satisfied; or
$$
S^j \text{ intersects both } CO^2(E^s_3) \text{ and } CO^2(-E^s_3), \tag 30
$$
when (28) holds. Then the system $\Sigma^{pr}$ is controllable.
\endproclaim

\medskip

The most simple case is when we have two systems $\Sigma_j$ and $Sigma_l$ with complex eigenvalues with different types of dynamics: when
$\Sigma^{pr}_j$ has the attracting cycle $S^j$, while $\Sigma^{pr}_l$ has the repulsive cycle $S^l$. Then, using the system $\Sigma^{pr}_j$ with
attracting cycle we may construct a trajectory from an arbitrary point $q$ to an arbitrary point of the attractive cycle $S^j$ - say, to its
intersection with the repulsive cycle $S^l$, and then continue this trajectory further to an arbitrary point with the help of the system
$\Sigma^{pr}_l$. We have the following controllability condition.

\medskip
\proclaim{Theorem~C} Let $\Sigma^{pr}$ be some projected system which satisfies the $LARC$ (CC1), the openness condition CC2) and generality
conditions CC3, CC4; and in the finite set ${\Cal U}^N=\{u_k|k=1,...,N\}$ of parameters for one of them $u_j$ the system $\Sigma_j$ has complex
eigenvalues satisfying (27), and for another $u_l$ the system $\Sigma_l$ has complex eigenvalues satisfying (28). Then the system $\Sigma^{pr}$ is
controllable.
\endproclaim

\enddocument
\bye